\documentclass[11pt]{amsart}
\usepackage{amssymb}
\usepackage{graphicx}

\def\N{\mathbb{N}}
\def\R{\mathbb{R}}
\def\C{\mathbb{C}}

\def\S{\mathbb{S}}
\def\disk{\mathbb{D}}
\def\Cc{\widehat{{\C}}}

\def\ol{\overline}
\def\tilde{\widetilde}

\def\and{and }

\def\crit{{\rm CV}}

\def\rho{\varrho}
\def\phi{\varphi}

\newtheorem{theorem}{Theorem}[section]
\newtheorem{proposition}[theorem]{Proposition}
\newtheorem{lemma}[theorem]{Lemma}
\newtheorem{definition}[theorem]{Definition}
\newtheorem{corollary}[theorem]{Corollary}
\newcounter{claiming}
\newenvironment{claim}{\refstepcounter{claiming}\par\medskip\noindent\textit{Claim \arabic{claiming}:}}{}

\title[Postcritically Fixed Newton Maps]{A Combinatorial Classification of Postcritically Fixed Newton Maps}
\author{Johannes R\"uckert}
\thanks{The author was partly supported by a Doktorandenstipendium of the German Academic Exchange Service (DAAD)}
\address{International University Bremen, School of Engineering \& Science, Campus Ring 12, 28759 Bremen, Germany}
\email{j.rueckert@iu-bremen.de}

\begin{document}
\begin{abstract}
We give a combinatorial classification for the class of postcritically fixed Newton maps of polynomials and indicate potential for extensions. 

As our main tool, we show that for a large class of Newton maps that includes all hyperbolic ones, every component of the basin of an attracting fixed point can be connected to $\infty$ through a finite chain of such components.
\end{abstract}
\maketitle

%%%%%%%%%%%%%%%%%%%%%%%%%%%%%%%%%%%%%%%
\section{Introduction} 

One of the most important open problems in rational dynamics is understanding the structure of the space of rational functions of a fixed degree $d\geq 2$. This problem is today wide open. 

Aside from being a useful tool for numerical root-finding, Newton maps of polynomials form an interesting subset of the space of rational maps that is more accessible for studying than the full space of rational maps. Hence, a partial goal in the classification of {\em all} rational maps can be to gain an understanding of the space of Newton maps. 

In this paper, we present a theorem that structures the dynamical plane of postcritically finite Newton maps, and then use this result to construct a graph that classifies those Newton maps whose critical orbits all terminate at fixed points. Newton maps of degree $1$ and $2$ are trivial, and we exclude these cases from our investigation.

\begin{definition}[Immediate Basin]
\label{Def_ImmediateBasin}
Let $f$ be a Newton map and $\xi\in\C$ a fixed point of $f$. Let $B_{\xi}=\{z\in\C\,:\, \lim_{n\to\infty}f^{\circ n}(z)=\xi\}$ be the {\em basin (of attraction)} of $\xi$.
The component of $B_\xi$ containing $\xi$ is called the {\em immediate basin} of $\xi$ and denoted $U_\xi$.
\end{definition}
Clearly, $B_\xi$ is open and by a theorem of Przytycki \cite{Przytycki}, $U_{\xi}$ is simply connected and unbounded (in fact, a result of Shishikura \cite{Shishikura} implies that every component of the Fatou set is simply connected). Moreover, $\infty\in\partial U_\xi$ is an accessible boundary point.

Our first result is the following.
\begin{theorem}[Preimages Connected]
\label{Thm_PreimagesConnected}
Let $f:\Cc\to\Cc$ be a Newton map with attracting fixed points $\xi_1,\dots,\xi_d\in\C$, and let $U'_0$ be a component of some $B_{\xi_i}$. Then, $U'_0$ can be connected to $\infty$ by the closures of finitely many components $U'_1,\dots, U'_k$ of $\bigcup_{i=1}^d B_{\xi_i}$. 

More precisely, there exists a curve $\gamma:[0,1]\to\Cc$ such that $\gamma(0)=\infty$, $\gamma(1)\in U'_0$ and for every $t\in[0,1]$, there exists $m\in\{0,1,\dots, k\}$ such that $\gamma(t)\in \ol{U}'_m$.
\end{theorem}

We will see that $\gamma$ can be chosen to consist of the closures of internal rays in the $U'_m$.
Theorem \ref{Thm_PreimagesConnected} allows to describe how the components of the basins are connected to each other. Thus, it is a basis for a combinatorial classification of certain Newton maps: Theorems \ref{Thm_NewtonGraph} and \ref{Thm_Realization} show that the combinatorics of these connections suffice to describe postcritically fixed Newton maps uniquely.

We call a Newton map {\em postcritically fixed} if all its critical points are mapped onto fixed points after finitely many iterations. If $f$ is a postcritically fixed Newton map, Theorem \ref{Thm_PreimagesConnected} allows to structure the entire Fatou set, because each Fatou component is in the basin of some attracting fixed point. Then, we construct the channel diagram $\Delta$ of $f$ (see Section \ref{Sec_Models}) and pull it back several times to get a connected graph $\Gamma$ that contains the forward orbits of all critical points, similar to the Hubbard tree of a postcritically finite polynomial. 

Conversely, we show that for every {\em abstract Newton graph} (a graph with dynamics that satisfies several natural conditions, see Definition \ref{Def_NewtonGraph}), there exists a unique postcritically fixed Newton map realizing it. The assignments of a Newton map to an abstract Newton graph and vice versa are injective and inverse to each other, so we give a combinatorial classification of postcritically fixed Newton maps by way of abstract Newton graphs. Thus, our main results are the following (see Sections \ref{Sec_Graph} and \ref{Sec_Thurston} for the precise definitions).

\begin{theorem}[Newton Map Generates Newton Graph]
\label{Thm_NewtonGraph}
Every postcritically fixed Newton map $f$ gives rise to a unique abstract Newton graph.
More precisely, there exists a unique $N\in\N$ such that $(\Delta_N,f)$ is an abstract Newton graph.

If $f_1$ and $f_2$ are Newton maps with channel diagrams $\Delta_1$ and $\Delta_2$ such that $(\Delta_{1,N},f_1)$ and $(\Delta_{2,N},f_2)$ are equivalent as abstract Newton graphs, then $f_1$ and $f_2$ are affinely conjugate. 
\end{theorem}
\begin{theorem}[Newton Graph Generates Newton Map]
\label{Thm_Realization}
Every abstract Newton graph is realized by a postcritically fixed Newton map which is unique up to affine conjugacy.
More precisely, let $(\Gamma,g)$ be an abstract Newton graph. Then, there exists a postcritically fixed Newton map $f$ with channel diagram $\tilde{\Delta}$ such that $(\ol{g},\Gamma')$ and $(f,\tilde{\Delta}'_{N_{\Gamma}})$ are Thurston equivalent as marked branched coverings. 

Moreover, if $f$ realizes two abstract Newton graphs $(\Gamma_1,g_1)$ and $(\Gamma_2,g_2)$, then the two abstract Newton graphs are equivalent.
\end{theorem}

Our construction of an abstract Newton graph can be done for all postcritically finite Newton maps, but will in general not contain the orbits of {\em all} critical points, and thus not describe the combinatorics of the entire Fatou set (note that there are rational maps with {\em buried} Fatou components that are not attached to any other Fatou component (it is not hard to find Newton maps with this property either). An extreme example of this behavior is provided by rational maps with Sierpinski Julia sets, see e.g.\ \cite[Appendix F]{Milnor2}). It seems likely however that with Theorem \ref{Thm_PreimagesConnected} and additional combinatorial objects that describe any strictly periodic or preperiodic critical points, a classification of at least all hyperbolic Newton maps can be achieved. Thus, our results are a first step towards a combinatorial classification of Newton maps, and in particular of all hyperbolic components in the space of Newton maps. They may also be a basis for transporting the powerful concept of Yoccoz puzzles, which has been used to prove local connectivity of the Julia set for many classes of polynomials, to the setting of Newton maps beyond the cubic case (Roesch has successfully applied Yoccoz puzzles to cubic Newton maps \cite{Roesch}).

A number of people have studied Newton maps and used combinatorial models to structure the parameter spaces of some Newton maps. Janet Head \cite{Head} introduced the {\em Newton tree} to characterize postcritically finite cubic Newton maps. Tan Lei \cite{TanLei}  built upon this work and gave a classification of postcritically finite cubic Newton maps in terms of matings and captures.  Jiaqi Luo \cite{Luo} extended some of these results to ``unicritical'' Newton maps, i.e.~Newton maps of arbitrary degree with only one {\em free} (non-fixed) critical value. 
The present work can be seen as an extension of these results beyond the setting of a single free critical value. The main differences to this setting are that the channel diagram is in general not a tree anymore and that in the presence of more than one non-fixed critical value, the iterated preimages of the channel diagram may be disconnected. 

\smallskip
This article is structured as follows. In Section \ref{Sec_Models}, we introduce the concept of a channel diagram for Newton maps and discuss some of its properties. We use the channel diagram and its preimages to prove Theorem \ref{Thm_PreimagesConnected} in Section \ref{Sec_Proof}. In Section \ref{Sec_Graph}, we introduce abstract Newton graphs and prove Theorem \ref{Thm_NewtonGraph}. Theorem \ref{Thm_Realization} is proved in Section \ref{Sec_Thurston}, following a review of some aspects of Thurston theory. We also give an introduction to the combinatorics of {\em arc systems} and state a result by Kevin Pilgrim and Tan Lei that restricts the possibilities of how arc systems and Thurston obstructions can intersect.

\subsection{Notation} 
Let us make precise what we mean by a Newton map.
\begin{definition}[Newton Map]
\label{Def_NewtonMap}
A rational function $f:\Cc\to\Cc$ of degree $d\geq 3$ is called a {\em Newton map} if $\infty$ is a repelling fixed point of $f$ and for each fixed point $\xi\in\C$, there exists an integer $m\geq 1$ such that $f'(\xi)=(m-1)/m$. 
\end{definition}

This definition is motivated by the following observation, which is a special case of \cite[Proposition 2.8]{RS} (the case of superattracting fixed points, i.e.\ every $m=1$, goes back to \cite[Proposition 2.1.2]{Head}).
\begin{proposition}[Head's Theorem]
\label{Prop_Head}
A rational map $f$ of degree $d\geq 3$ is a Newton map if and only if there exists a polynomial $p:\C\to\C$ such that for $z\in\C$, $f(z)=z-p(z)/p'(z)$. 
\qed 
\end{proposition}
Let $f$ be a Newton map. A point $z\in\C$ is called a {\em pole} if $f(z)=\infty$ and a {\em prepole} if $f^{\circ k}(z)=\infty$ for some minimal $k>1$.
If $g:\S^2\to\S^2$ is a branched covering map, we call a point $z\in\S^2$ a {\em critical point} if $g$ is not injective in any neighborhood of $z$. For the Newton map $f$, this is equivalent to saying that $z\in\C$ and $f'(z)=0$, because $\infty$ is never a critical point of $f$. It follows from the Riemann-Hurwitz formula \cite[Theorem 7.2]{Milnor} that a degree-$d$ branched covering map of $\S^2$ has exactly $2d-2$ critical points, counting multiplicities.
\begin{definition}[Postcritically Fixed]
\label{Def_PostCritFixed}
Let $g:\S^2\to\S^2$ be a branched covering map of degree $d\geq 2$ with (not necessarily distinct) critical points $c_1,\dots,c_{2d-2}$. We denote the set of {\em critical values} of $g$ by 
\[
\crit(g) :=\{g(c_1),\dots,g(c_{2d-2})\}\;.
\]
Then, $g$ is called {\em postcritically finite} if the set
\[
	P_g := \bigcup_{n\geq 0} g^{\circ n}(\crit(g))
\]
is finite. We say that $g$ is {\em postcritically fixed} if there exists $N\in\N$ such that for each $i\in\{1,\dots,2d-2\}$, $g^{\circ N}(c_i)$ is a fixed point of $g$.
\end{definition}

\begin{definition}[Access to $\infty$]
\label{Def_Access}
Let $V\subset \Cc$ be a simply connected unbounded domain and consider a curve $\Gamma:[0,\infty)\to V\cup\{\infty\}$ with $\Gamma(0)\in V$ and $\lim_{t\to\infty}\Gamma(t)=\infty$, such that $\Gamma(t)\neq\infty$ for all $t<\infty$. Its homotopy class within $V$ defines an {\em access to $\infty$} for $V$. In other words, a curve $\Gamma'$ with the same properties lies in the same access as $\Gamma$ if the two curves are homotopic in $V$, fixing the endpoint $\infty$.

If $V=U_{\xi}$ is an immediate basin, then we require that $\Gamma(0)=\xi$ and $\Gamma$ and $\Gamma'$ are homotopic with both endpoints fixed.
\end{definition}
\begin{proposition}[Accesses]
\label{Prop_Access}
({\em c.f.\ \cite{HSS}}) Let $f$ be a Newton map of degree $d\geq 3$ and $U_\xi$ an immediate basin for $f$. Then, there exists $k\in \{1,\dots,d-1\}$ such that $U_\xi$ contains $k$ critical points of $f$ (counting multiplicities), $f|_{U_\xi}$ is a covering map of degree $k+1$, and $U_\xi$ has exactly $k$ accesses to $\infty$.
\qed
\end{proposition}

%%%%%%%%%%%%%%%%%%%%%%%%%%%%
\section{The Channel Diagram}
\label{Sec_Models}
In the following, by a {\em (finite) graph} we mean a connected topological space $\Gamma$ homeomorphic to the quotient of a finite disjoint union of closed arcs by an equivalence relation on the set of their endpoints. The arcs are called {\em edges} of the graph, an equivalence class of endpoints a {\em vertex}. We usually consider imbedded graphs in $\S^2$, i.e.\ the homeomorphic image of a graph in $\S^2$. 

\begin{definition}[Graph Map]
Let $\Gamma_1,\Gamma_2\subset\S^2$ be two finite graphs and $g:\Gamma_1\to\Gamma_2$ continuous. We call $g$ a {\em graph map} if it is injective on each edge of $\Gamma_1$ and forward and inverse images of vertices are vertices. If the graph map $g$ is a homeomorphism, then we call it a {\em graph homeomorphism}.
\end{definition}

In the following, the closure and boundary operators will be understood with respect to the topology of $\Cc$, unless otherwise stated. Also, we will say that a set $X\subset\Cc$ is {\em bounded} if $\infty\not\in\ol{X}$.

We say that a Newton map $f$ of degree $d\geq 3$ \emph{satisfies \eqref{Prop_PCF}} if it has the following property:
\begin{equation}
\tag{$\star$}
\label{Prop_PCF}
\left\{
\begin{aligned}
& \text{if } c \text{ is a critical point of } f \text{ with } c\in B_{\xi_j} \text{ for some } \\
& j\in\{1,\dots,d\}, \text{ then } c \text{ has finite orbit.}
 \end{aligned}
 \right.
\end{equation}

We omit the easy proof of the following well-known fact.
\begin{lemma}[Only Critical Point]
\label{Lem_OnlyCritical}
Let $f$ be a Newton map that satisfies \eqref{Prop_PCF} and let $\xi\in\C$ be a fixed point of $f$ with immediate basin $U_{\xi}$. Then, $\xi$ is the only critical point in $U_{\xi}$.
\qed 
\end{lemma}

\begin{figure}[hbt]
\begin{center}
\setlength{\unitlength}{1cm}
\begin{picture}(10,10)
\includegraphics[viewport=0 0 412 412,width=10cm]{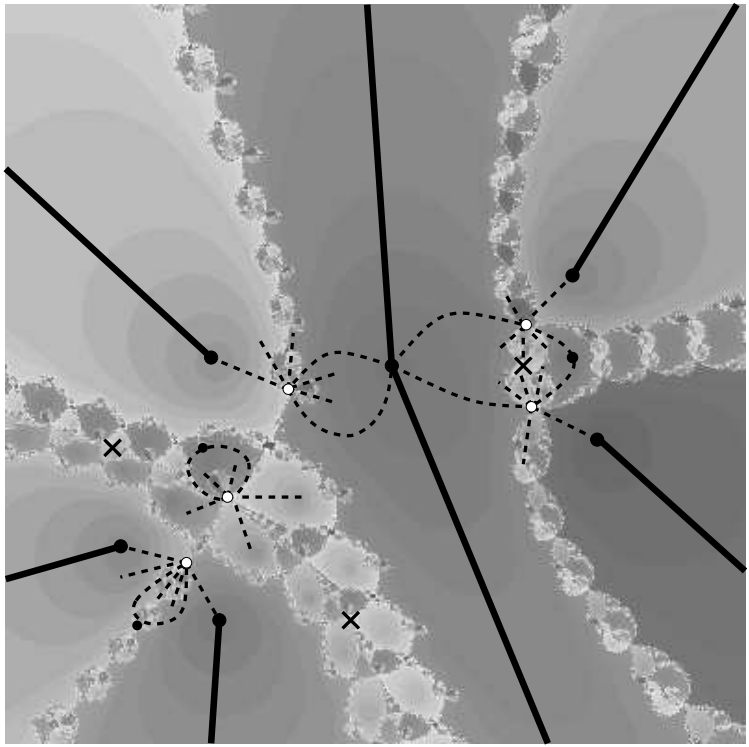}
\end{picture}
\caption{\label{Fig_Degree6} A Newton map of degree 6, superimposed with its channel diagram: the solid lines represent the fixed rays of the immediate basins, the black dots correspond to the fixed points (the vertex at $\infty$ is not visible). 
The dashed lines show the first preimage of the channel diagram: white circles represent poles, a cross is a free critical point. The Newton map has a pole outside $\Delta_1$, and the right boundary component of the central immediate basin contains two poles.}
\end{center}
\end{figure}

It follows that if $f$ is a Newton map that satisfies \eqref{Prop_PCF}, then each immediate basin $U_{\xi}$ has a global B\"ottcher map $\phi_\xi: (\disk,0)\to (U_\xi,\xi)$ with the property that $f(\phi_\xi(z)) = \phi_\xi(z^{k_\xi})$ for each $z\in \disk$, where $k_\xi-1\geq 1$ is the multiplicity of $\xi$ as a critical point of $f$ \cite[Theorems 9.1 \& 9.3]{Milnor}. The $k_\xi-1$ radial lines (or {\em internal rays}) in $\disk$ which are fixed under $z\mapsto z^{k_\xi}$ map under $\phi$ to $k_\xi-1$ pairwise disjoint, non-homotopic injective curves $\Gamma_{\xi}^1,\dots,\Gamma_{\xi}^{k_\xi}$ in $U_{\xi}$ that connect $\xi$ to $\infty$ and are each invariant under $f$. They represent all accesses to $\infty$ of $U_{\xi}$, see Proposition \ref{Prop_Access}. Hence if $\xi_1,\dots,\xi_d\in\C$ are the attracting fixed points of $f$, then the union 
\[
	\Delta := \bigcup_{i=1}^d\bigcup_{j=1}^{k_{\xi_i}} \ol{\Gamma_{\xi_i}^j}
\]
of these invariant curves over all immediate basins forms a connected and $f$-invariant graph in $\Cc$ with vertices at the $\xi_i$ and at $\infty$.  We call $\Delta$ the {\em channel diagram} of $f$.
The channel diagram records the mutual locations of the immediate basins of $f$ and provides a first-level combinatorial structure to the dynamical plane. Figure \ref{Fig_Degree6} shows a Newton map and its channel diagram. The following definition is an axiomatization of the channel diagram.

\begin{definition}[Abstract Channel Diagram]
\label{Def_ChannelDiagram} An {\em abstract channel diagram} of degree $d\geq 3$
is a graph $\Delta \subset \S^2$ with vertices $v_0,\dots,v_d$ and
edges $e_1,\dots,e_l$ that satisfies the following properties:
\begin{enumerate}
\item $l\leq 2d-2$;
\item each edge joins $v_0$ to a $v_i$, $i>0$;
\item each $v_i$ is connected to $v_0$ by at least one edge;
\item \label{necessaryCondition} if $e_i$ and $e_j$ both join $v_0$ to $v_k$, then each connected component of
$\S^2\setminus \ol{e_i\cup e_j}$ contains at least one vertex of $\Delta$.
\end{enumerate}
\end{definition}

We say that an abstract channel diagram $\Delta$ is {\em realized} if there exist
a Newton map with channel diagram $\hat{\Delta}$ and a graph homeomorphism $h:\Delta\to\hat{\Delta}$ that preserves the cyclic order of edges at each vertex.

We claim that if $\hat{\Delta}$ is the channel diagram of the Newton map $f$, then it is an abstract channel diagram: by construction, $\hat{\Delta}$ has at most $2d-2$ edges and it satisfies (2) and (3). Finally, $\hat{\Delta}$ satisfies (\ref{necessaryCondition}), because for any immediate basin $U_{\xi}$ of $f$, every component of $\C\setminus U_\xi$ contains at least one fixed point of $f$ \cite[Corollary 5.2]{RS} (see also Theorem \ref{Thm_FixedPoles} below).

The following useful observation is a special case of \cite[Theorem 4.8]{RS}.
\begin{lemma}[Fixed Points]
\label{Lem_Lefschetz}
Let $f$ be a Newton map and let $D\subset\Cc$ be a closed topological disk such that $\gamma:=f(\partial D)$ is a simple closed curve with the property that $\gamma\cap\mathring{D}=\emptyset$. Let $V$ be the unique component of $\Cc\setminus\gamma$ that contains $\mathring{D}$ and let $\{\gamma'_i\}_{i\in I}$ be the collection of boundary components of $f^{-1}(V)\cap D$. Then, the number of fixed points of $f$ in $D$ equals
\[
	\sum_{i\in I} \left|\deg(f|_{\gamma'_i}:\gamma'_i\to\gamma)\right|\;.
\]
In particular, if $f^{-1}(V)\cap D\neq\emptyset$, then $D$ contains a fixed point.
\qed
\end{lemma}
\begin{remark}
Since $f$ has no parabolic fixed points, we do not need to take multiplicities of fixed points into account.

Note also that the $\gamma'_i$ are exactly the components of $f^{-1}(\gamma)\cap D$, except possibly $\partial D$ itself. The boundary is excluded if points in $\mathring{D}$ near $\partial D$ are mapped out of $\ol{V}$. By the lemma, the only case in which $D$ does not contain a fixed point of $f$ is if all of $D$ is mapped outside of $\ol{V}$. 
\end{remark}

The following theorem shows a relation between poles and fixed points outside immediate basins. It considerably sharpens \cite[Corollary 5.2]{RS}, which states that for an immediate basin $U_{\xi}$ of a Newton map, every component of $\C\setminus U_{\xi}$ cointains at least one fixed point. 

\begin{theorem}[Fixed Points and Poles]
\label{Thm_FixedPoles}
Let $f$ be a Newton map and $U_\xi$ an immediate basin. If $W$ is a component of $\C\setminus U_\xi$, then the number of fixed points in $W$ equals the number of poles in $W$, counting multiplicities. 
\end{theorem}
\begin{proof} 
Let $d\geq 3$ be the degree of $f$. If $U_\xi$ does not separate the plane, i.e.\ if it has only one access to $\infty$, then the claim follows trivially: $W$ contains all $d-1$ finite poles and the $d-1$ other finite fixed points of $f$. So suppose in the following that there is a Riemann map $\phi:(\disk,0)\to (U_\xi,\xi)$ with $f(\phi(z))=\phi(g(z))$ for some Blaschke product $g:\disk\to\disk$ of degree $k:=\deg(f|_{U_{\xi}}:U_{\xi}\to U_{\xi})\geq 3$. 
\begin{figure}[hbt]
\begin{center}
\setlength{\unitlength}{1cm}
\begin{picture}(9,5.2)
\put(0,0){\includegraphics[viewport=000 000 413 238,width=9cm]{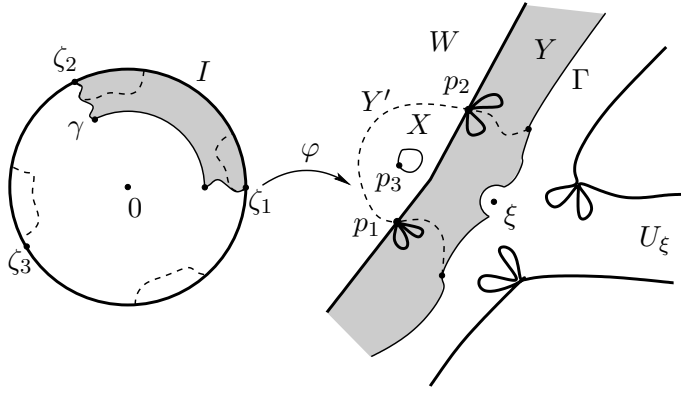}}
\put(3.9,3.1){$\phi$}
\put(1.6,2.3){$0$}
\put(0.8,3.4){$\gamma$}
\put(8.4,1.9){$U_{\xi}$}
\put(7.5,4){$\Gamma$}
\put(6.6,2.2){$\xi$}
\put(3.2,2.4){$\zeta_1$}
\put(0,1.6){$\zeta_3$}
\put(0.6,4.3){$\zeta_2$}
\put(5.6,4.5){$W$}
\put(2.5,4.1){$I$}
\put(4.6,2.1){$p_1$}
\put(5.8,4){$p_2$}
\put(7,4.4){$Y$}
\put(4.7,3.7){$Y'$}
\put(5.3,3.4){$X$}
\put(4.9,2.7){$p_3$}
\end{picture}
\caption{\label{Fig_LefschetzCurve} In the proof of Theorem \ref{Thm_FixedPoles}, the construction of $\gamma\subset\disk$ is shown on the left for $k=4$. The dashed curves are the components of 
$g^{-1}(\gamma)$. The right picture shows the curve $\Gamma\subset U_{\xi}$. The dashed curve indicates where $\Gamma'_1$ differs from $\Gamma$.}
\end{center}
\end{figure}

We may extend $g$ by reflection to a rational function $\hat{g}:\Cc\to\Cc$ of degree $k$ whose Julia set equals $\S^1$ and that has $k-1$ fixed points $\zeta_1,\dots,\zeta_{k-1}\in \S^1$. These fixed points correspond to the accesses to $\infty$ of $U_{\xi}$. Since $\hat{g}$ fixes $\disk$, the $\zeta_i$ have real positive multipliers and since $0$ and $\infty$ attract all of $\disk$ and of 
$\Cc\setminus\ol{\disk}$, respectively, none of the $\zeta_i$ can be attracting or parabolic. Hence they are pairwise distinct and repelling. For each $\zeta_i$, choose a linearizing neighborhood and choose $0<\rho<1$ large enough so that all critical values of $\hat{g}$ in $\disk$ have absolute value less than $\rho$ and so that the linearizing neighborhoods of all $\zeta_i$ intersect the circle at radius $\rho$. Let $\gamma\subset\disk$ be the unique curve with the following properties, see Figure \ref{Fig_LefschetzCurve}: there are adjacent fixed points $\zeta_j,\zeta_{j+1}\in\S^1$ and injective curves $\gamma_j$, $\gamma_{j+1}$, so that $\gamma_j$ connects $\zeta_j$ to the circle at radius $\rho$ and $\gamma_j$ is a straight line segment in linearizing coordinates of $\zeta_j$; the same for $\gamma_{j+1}$ and $\zeta_{j+1}$. Their closures separate the circle at radius $\rho$ into two arcs. Of those arcs, let $\gamma'$ be the one for which $\gamma:=\gamma_j\cup \ol{\gamma'}\cup \gamma_{j+1}$ has the property that $\phi(\gamma)\subset U_{\xi}$ separates $W$ from $\xi $. Let $\Gamma:=\phi(\gamma)\cup\{\infty\}$. Then, $\Gamma$ is a simple closed curve in $\ol{U}_{\xi}$ and contains no critical values, except possibly $\infty$. 

Let us first suppose that $\infty$ is not a critical value. Then, every component $\Gamma'_i$ of $f^{-1}(\Gamma)$ is a simple closed curve and $\deg(f|_{\Gamma'_i}:\Gamma'_i\to\Gamma)$ equals the number of poles on $\Gamma'_i$ (here, we do not need to count multiplicities, because we have assumed that $f$ has no critical poles). 

Let $\Gamma_1'$ be the component of $f^{-1}(\Gamma)$ containing $\infty$. We claim that $\Gamma'_1\cap U_{\xi}$ consists of two connected components, each of which is an injective curve that connects $\infty$ to a pole on $\partial U_{\xi}$; call these poles $p_1$ and $p_2$. Indeed, consider the situation in $\disk$-coordinates. Let $I\subset\S^1$ be the arc between $\zeta_j$ and $\zeta_{j+1}$ that is separated from $0$ by $\gamma$. Since $\mathring{I}$ contains no fixed points of $\hat{g}$, $\hat{g}(\ol{I})$ covers $\S^1\setminus \ol{I}$ exactly once and $\ol{I}$ itself exactly twice. Hence, it is easy to see that $\hat{g}^{-1}(\gamma)$ has exactly two connected components that intersect $\gamma$. This proves the claim. 

Let $Y$ be the closure of the component of $\Cc\setminus \Gamma$ that contains $W$ and let $Y'$ be the closure of the component of $\Cc\setminus\Gamma'_1$ that intersects $W$ in an unbounded set. We distinguish two cases.

If $p_1=p_2$, then $\Gamma'_1\subset \ol{U}_{\xi}$ and $W\subset Y'$. Moreover, $Y'$ is a closed topological disk that contains exactly the same fixed points and poles of $f$ as $W$. Since $f(\partial Y')\cap \mathring{Y}'=\emptyset$, Lemma \ref{Lem_Lefschetz} gives that the number of fixed points in $Y'$ (including $\infty$) equals the number of poles (again including $\infty$), because on every component of $f^{-1}(\Gamma)$ in $Y'$, the degree of $f$ equals the number of poles it contains. Excluding $\infty$ again, the claim follows.

If $p_1\neq p_2$, then $\Gamma'_1\not\subset \ol{U}_{\xi}$ and $W\not\subset Y'$ (this is the situation pictured in Figure \ref{Fig_LefschetzCurve}). If the set $X:=W\setminus Y'$ contains neither poles nor fixed points of $f$, then we can proceed as before. 

Indeed, $\Gamma'_1\cup U_{\xi}$ separates $X$ from $\infty$ and since every fixed point of $f$ is surrounded by its unbounded immediate basin, $X$ cannot contain a fixed point. Now suppose by way of contradiction that $X$ contains a pole $p_3$ of $f$. If $p_3\in \partial X$, then $p_3\in\partial U_{\xi}$. But this would imply the existence of an additional pre-fixed point of $\hat{g}$ on $I$, a contradiction. The other case is that $p_3\in\mathring{Y}''$. Then, there exists a component $\Gamma'_2$ of $f^{-1}(\Gamma)$ in $\mathring{Y}''$. Let $D$ be the bounded disk bounded by $\Gamma'_2$. We may assume without loss of generality that there is no component of $f^{-1}(\Gamma)$ separating $D$ from $\partial X$. Since points in $Y'$ near $\Gamma_1'$ are mapped into $Y$ under $f$, points in $X$ near $\Gamma_1'$ are mapped out of $Y$, and it follows that again $f(D)\cap Y\neq\emptyset$. Now, $\ol{D}$ contains a fixed point by Lemma \ref{Lem_Lefschetz}. This is a contradiction.

In the remaining case that $\infty$ is a critical value, we perturb $f$ slightly to avoid that situation. Since poles and fixed points of $f$ move continuously under perturbation, and $\Gamma$ does too, this does not change the count. Note that while $U_{\xi}$ and $W$ might move discontinuously, this does not pose a problem because we have actually counted poles and fixed points in $Y'$, after having established that we do not lose anything by this replacement.
\end{proof}
\begin{corollary}[Fixed Points in Complement]
\label{Cor_FixedComplement}
Let $f$ be a Newton map that satisfies \eqref{Prop_PCF} and let $\Delta$ be the channel diagram of $f$. Let $V$ be a component of $\Cc\setminus\Delta$ and let $p$ be the number of poles of $f$ in $V$, counting multiplicities. Then $\partial V \cap \C$ contains $p+1$ fixed points.
\end{corollary}
\begin{proof}
If $V$ is the only component of $\Cc\setminus\Delta$, the claim follows trivially. If $\xi$ is the only fixed point on $\partial V$ whose immediate basin $U_\xi$ separates the plane, then the claim follows directly from Theorem \ref{Thm_FixedPoles}. Indeed, let in this case $R_1,R_2$ be the fixed internal rays of $U_\xi$ that are on $\partial V$ and let $V_1$ be the component of $\C\setminus (R_1\cup R_2\cup\{\xi\})$ such that $V\subset V_1$. Then, $V_1$ also contains $p$ poles and by Theorem \ref{Thm_FixedPoles}, $V_1$ contains $p$ fixed points. Since $\xi\in\partial V$ as well, the claim follows.
\begin{figure}[hbt]
\begin{center}
\setlength{\unitlength}{1cm}
\begin{picture}(6,2.8)
\put(0,0){\includegraphics[viewport=0 0 337 158,width=6cm]{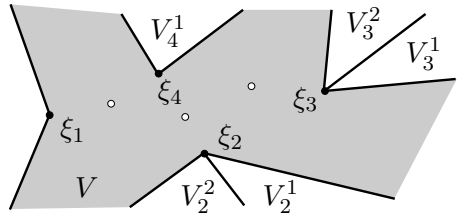}}
\put(0.9,0.2){$V$}
\put(0.7,1){$\xi_1$}
\put(2.8,0.9){$\xi_2$}
\put(3.8,1.4){$\xi_3$}
\put(2,1.5){$\xi_4$}
\put(1.9,2.3){$V_4^1$}
\put(5.3,2){$V_3^1$}
\put(4.5,2.4){$V_3^2$}
\put(3.4,0.1){$V_2^1$}
\put(2.3,0.1){$V_2^2$}
\end{picture}
\caption{\label{Fig_FixedComplement} In Corollary \ref{Cor_FixedComplement}, the open set $V$ is bounded by parts of $\Delta$. The black dots are fixed points, the white dots represent poles. The $V_j^i$ may well contain further structure of $\Delta$.}
\end{center}
\end{figure}

Now suppose that $\xi_1,\dots,\xi_k\in \partial V$ are the fixed points on $\partial V$ whose immediate basins separate the plane. Let $R_1, R_2$ be the fixed internal rays of $U_{\xi_1}$ on $\partial V$ and let $V_1$ be as above. Let $m$ be the number of poles in $V_1$. As before, it follows that $V_1$ contains $m$ fixed points. Let $m'=m-p$. For  $j=2,\dots,k$, denote by $V_j^1,\dots,V_j^{i_j}$ all complementary components of the closures of the fixed internal rays of $U_{\xi_j}$ that do not contain $V$. By Theorem \ref{Thm_FixedPoles}, each $V_j^i$ contains as many poles as fixed points, hence all $V_j^i$ combined contain $m'$ poles and $m'$ fixed points. Hence, $V_1\setminus(\bigcup_{j=2}^k\bigcup_{\ell=1}^{i_j} V_j^{\ell})$ contains $m-m'=p$ fixed points. The claim now follows, because including $\xi_1$, $\partial V\cap \C$ contains $p+1$ fixed points.
\end{proof}

\begin{corollary}[Existence of Shared Poles]
\label{Cor_SharedPoles}
Let $f$ be a Newton map that satisfies \eqref{Prop_PCF} and let $\Delta$ be the channel diagram of $f$. If $V$ is a component of $\Cc\setminus\Delta$, then there is at least one pair of fixed points $\xi_1,\xi_2\in\partial V\cap\C$ such that $\partial U_{\xi_1}$ and $\partial U_{\xi_2}$ intersect in a pole.
\end{corollary}
\begin{proof}
Let $U_{\xi}$ be an immediate basin. Clearly, the components of $\partial U_{\xi}\cap\C$ are separated by the accesses to $\infty$. We have seen that in the conjugate dynamics $\hat{g}|_{\disk}$, for every arc $I\subset\S^1$ between two fixed points, $\hat{g}(I)=\S^1$. Therefore, $I$ contains pre-fixed points of $\hat{g}$. Since poles and $\infty$ are accessible boundary points of $U_{\xi}$, we conclude that every component of $\partial U_{\xi}\cap\C$ contains at least one pole. By Corollary \ref{Cor_FixedComplement}, there has to be at least one pole in $V$ that is on the boundary of at least two immediate basins.
\end{proof}
\begin{remark}
 Figure \ref{Fig_Degree6} shows that a component of $\partial U_{\xi}\cap \C$ may contain more than one pole.

Note also that a simple pole is on the boundary of at most two immediate basins, because otherwise $f$ cannot preserve the cyclic order of the immediate basins near that pole. This was first observed by Janet Head \cite{Head}.
\end{remark}

%%%%%%%%%%%%%%%%%%%%%%%%%%%%
\section{Proof of Theorem \ref{Thm_PreimagesConnected}}
\label{Sec_Proof}

Let $\xi_1,\dots,\xi_d$ be the fixed points of $f$ and let $U'_0$ be a component of some $B_{\xi_i}$. 
For the main part of this section, we assume that $f$ satisfies the finiteness condition  \eqref{Prop_PCF}. We will indicate at the end of this section how to prove the result in the general case.

Let $\Delta$ be the  channel diagram of $f$. Recall that it consists of invariant rays within the $U_{\xi_i}$ that connect the $\xi_i$ to $\infty$. Denote by $\Delta_n$ the connected component of $f^{-n}(\Delta)$ that contains $\Delta$ (with this convention, $\Delta=\Delta_0$). Every edge of $\Delta_n$ is then an internal ray of a component of some $B_{\xi_i}$, while every vertex is a preimage of a $\xi_i$, or a pole or prepole.

To prove Theorem \ref{Thm_PreimagesConnected}, it suffices to show that there exists $n\in\N$ such that $\Delta_n$ contains all poles of $f$: then, every pole of $f$ can be connected to $\infty$ through a finite chain of internal rays in the $B_{\xi_i}$, and hence through a finite chain of components of the basins. By induction, each prepole is in $\Delta_m$ for sufficiently large $m$.
Since $\infty$ is on the boundary of every immediate basin, $\partial U'_0$ contains a prepole. This finishes the argument.

It remains to show that there exists $n\in\N$ such that $\Delta_n$ contains all poles of $f$.
If $\Delta_1$ contains all poles of $f$, then we are done. So assume in the following that there exists a component $C_1$ of $f^{-1}(\Delta)$ such that $C_1\cap\Delta_1=\emptyset$ (Figure \ref{Fig_Degree6} shows that this does occur). Equivalently, we may assume that there exists a component $V_0$ of $\Cc\setminus\Delta$ and a component $V_1$ of $f^{-1}(V_0)$ such that $V_1$ is multiply connected. Then, we choose $C_1$ so that it intersects $\partial V_1$.

Denote by $C_n$ the component of $f^{-n}(\Delta)$ containing $C_1$. We will assume that $C_n\cap\Delta_n=\emptyset$ for all $n\in\N$ (otherwise we would be done). We will lead this assumption to a contradiction.

\begin{lemma}[Preimage Inside]
\label{Lem_Inside}
With the above notation, $V_1\subset V_0$.
\end{lemma}
\begin{proof}
Since $\Delta\subset f^{-1}(\Delta)$, we either have $V_1\subset V_0$ or $V_1\cap V_0=\emptyset$. In the latter case, let $\gamma\subset V_0$ be a simple closed curve near $\partial V_0$ that surrounds all critical values within $V_0$ (note that $V_0$  contains critical values, because it is the image of the multiply connected domain $V_1$). Then, $f^{-1}(\gamma)\cap V_1$ consists of several nested and non-contractible (in $V_1$) simple closed curves. Let $\gamma'$ be the outermost of them and let $D$ be the bounded component of $\Cc\setminus \gamma'$. By Lemma \ref{Lem_Lefschetz}, $\ol{D}$ contains a fixed point of $f$. This is a contradiction, because $\ol{D}$ is separated from $\Delta$ and all fixed points are contained in $\Delta$.
\end{proof}

\begin{lemma}[Preimage Unbounded]
\label{Lem_V1}
Let $n\in\N$ and suppose that $W$ is an unbounded component of $\Cc\setminus\Delta_n$. If $W'$ is a component of $f^{-1}(W)$ with $W'\subset W$, then $W'$ is unbounded.
\end{lemma}
\begin{proof}
Let $\xi_1,\dots,\xi_k$ be the attracting fixed points of $f$ in $\partial W$. By Lemma \ref{Lem_OnlyCritical}, $f$ has B\"ottcher coordinates near each $\xi_i$. For $i\in\{1,\dots, k\}$, choose a neighborhood $B_i$ of $\xi_i$ that has the following properties: $B_i$ contains no critical values except $\xi_i$; in B\"ottcher coordinates, $B_i$ is a round disk centered at $0$, small enough so that $\Delta_{n+1}\cap B_i$ consists of radial lines; $f$ is conjugate to $z\mapsto z^{k_i}$ on $B_i$. Since $\Delta_n$ is a graph, there exists a simple closed curve $\gamma\subset \partial W\subset \Delta_n$ that surrounds $W$ and thus $W'$. By possibly modifying $\gamma$ within the $B_i$, we may assume that in B\"ottcher coordinates, every point $z\in \gamma\cap B_i$ is either in $\Delta_n$ or on a circle of constant radius centered at $0$. Let $D$ be the component of $\Cc\setminus \gamma$ that intersects $W$. Then, $D\cap \gamma \subset \bigcup_{i=1}^k B_i$.

Now suppose by way of contradiction that $W'\subset W$ is bounded and let $\gamma'$ be the outermost (in $\C$) simple closed curve in $f^{-1}(\gamma)$ that intersects $\partial W'$. Let $D'$ be the component of $\C\setminus\gamma'$ that is contained in $D$. Observe that $D'\cap\gamma=\emptyset$: $\gamma$ intersects $D'$ at most in the $B_i$. But there, $\gamma$ was chosen in such a way that $\gamma'$ is strictly further away from $\xi$ than $\gamma$, so in particular, $\gamma$ and $D'$ are disjoint. Now, either $f(D')\subset D$ or $\ol{D}'$ contains another component of $f^{-1}(\gamma)$. In both cases, Lemma \ref{Lem_Lefschetz} shows that $\ol{D}'$ contains a fixed point of $f$. This is a contradiction, because $\gamma$ was constructed in such a way that $\ol{D}'$ does not contain an attracting fixed point, while $\infty$ is not in $\ol{D}'$ by assumption.
\end{proof}

\begin{corollary}[Free Pole in Unbounded Nest]
\label{Cor_FreeNest}
Suppose that $\Delta_n\cap C_n=\emptyset$ for all $n\in\N$. Then for each $n\geq 1$, there exists an unbounded and multiply connected component $V_n$ of $f^{-n}(V_0)$ such that $V_n\subset V_{n-1}$ and $C_n\cap\partial V_n\neq\emptyset$, while the component of $\partial V_n$ that contains $\infty$ is in $\Delta_n$. 
\end{corollary}
\begin{proof}
By definition, $V_1$ is multiply connected and $\partial V_1\cap C_1\neq\emptyset$. By Lemma \ref{Lem_Inside}, $V_1\subset V_0$ and  by Lemma \ref{Lem_V1}, $V_1$ is unbounded.

Now suppose by induction that $V_{n}$ has the claimed properties. Clearly, every component of $f^{-1}(V_n)$ is multiply connected. Since $\infty \in \partial V_n$, there exists a component $V_{n+1}$ of $f^{-1}(V_n)$ such that $\partial V_{n+1}\cap C_{n+1}\neq\emptyset$. Then, $V_{n+1}\subset V_n$. 

Let $W$ be the component of $\Cc\setminus \Delta_n$ containing $V_n$ and let $W'$ be the component of $f^{-1}(W)$ containing $V_{n+1}$. If $V_{n+1}$ was bounded, then so would be $W'$, contradicting Lemma \ref{Lem_V1}. Since $\Delta_{n+1}$ is the only unbounded component of $f^{-(n+1)}(\Delta)$, we are done.
\end{proof}

We call the unbounded component of $\partial V_n$ the {\em outer boundary} and denote it with $B_n$.

Recall that to finish the proof of Theorem \ref{Thm_PreimagesConnected} in case of the finiteness condition \eqref{Prop_PCF}, it suffices to show the following.

\begin{theorem}[Poles Connect to $\infty$]
\label{Thm_PolesConnect}
There exists $n\in\N$ such that $\Delta_n$ contains all poles of $f$.
\end{theorem}
\begin{proof} Suppose by way of contradiction that $C_n\cap \Delta_n=\emptyset$ for all $n\in\N$. Suppose first that for large enough $n$, $V_n$ surrounds only one pole, and this pole is in $C_n$.

For any $n\in\N$, we call a closed arc in $B_n$ a \emph{bridge} if it connects two distinct finite fixed points $\xi_1$ and $\xi_2$ and contains no other fixed points, in particular not $\infty$. The \emph{length} of a bridge is the number of edges of $\Delta_n$ it consists of. We say that two distinct finite fixed points $\xi_1, \xi_2\in B_n$  are \emph{adjacent} if there exists a bridge in $B_n$ connecting them.

By Corollary \ref{Cor_FreeNest}, all $V_n$ are unbounded and hence have at least one access to $\infty$. Moreover, every unbounded edge of $B_n$ is contained in $\Delta$ and connects $\infty$ to a finite fixed point. It now follows from Definition \ref{Def_ChannelDiagram} \eqref{necessaryCondition}, that each $B_n$ contains at least two finite fixed points. Hence, the number of fixed points and accesses in $\partial V_n$ cannot shrink infinitely often and there exists a minimal $n_0$ such that for $n\geq n_0$, $B_n$ and $B_{n_0}$ contain the same fixed points and $V_n$ and $V_{n_0}$ have the same accesses to $\infty$. Therefore, two distinct finite fixed points $\xi_1,\xi_2\in B_n$ are adjacent for $n\geq n_0$ if and only if they are adjacent in $B_{n_0}$: a bridge in $B_n$ between them that does not exist in $B_{n_0}$ would separate another fixed point on $B_{n_0}$ (or an access to $\infty$ of $V_{n_0}$) from $V_n$; the other direction is trivial. 
\begin{claim}
\label{Claim_1}
If $\xi$ is a finite fixed point that is in $B_n$ for all $n\in\N$, then $\xi$ is adjacent to at most one finite fixed point.
\end{claim}
\begin{proof}[Proof of Claim \ref{Claim_1}:]
Let $\phi:\disk\to U_{\xi}$ be a B\"ottcher map that conjugates $f$ to $g(z)=z^k$ for some $k\geq 2$ and  let $R_1,\dots, R_{k-1}\subset\disk$ be the fixed internal rays of $g$. For all $n\in\N$, $\phi^{-1}(V_n\cap U_{\xi})$ equals a sector $S_n$ of $\disk\setminus \bigcup_{i=1}^{k-1} \ol{g^{-n}(R_i)}$, such that $S_{n+1}\subset S_n$ and $g(S_{n+1})= S_n$. It is easy to see that this can only happen if one of the boundary arcs of all $S_n$ is a fixed ray.

It follows that $\xi$ and $\infty$ are connected in all $B_n$ by a fixed edge, and at most one other edge that ends at $\xi$ can be part of a bridge in $B_n$.
\end{proof}
Observe also that $B_{n_0}$ contains a bridge: by minimality of $n_0$, there is a finite fixed point on $B_{n_0-1}$ or an access to $\infty$ of $V_{n_0-1}$ that is separated from $V_{n_0}$ by an arc in $B_{n_0}$. Possibly by extending this arc along $B_n$, we find a bridge $X_{n_0}\subset B_{n_0}$ between two finite fixed points (this is possible because the extension must hit a finite fixed point before $\infty$). If $n_0=1$, then $B_1$ contains a bridge $X_{n_0}$ between two finite fixed points by Corollary \ref{Cor_SharedPoles}. Let $\xi_1$ and $\xi_2$ be the endpoints of $X_{n_0}$. 

By induction on $n\geq n_0$, let $X_{n+1}$ be the shortest arc in $f^{-1}(X_n)$ that ends at $\xi_1$ and connects $\xi_1$ to a point $\xi'_2\in f^{-1}(\{\xi_2\})$ within $B_{n+1}$. If $\xi'_2=\xi_2$, then $X_{n+1}$ is a bridge as well, otherwise we say that the bridge $X_{n_0}$ \emph{breaks up} at time $n+1$. 
Note that if $X_{n+1}$ is a bridge, then the length of $X_{n+1}$ can only be greater than the one of $X_{n}$ if $X_{n+1}$ contains a pole. Hence the length is constant for all sufficiently large $n$. 

\begin{claim}
\label{Claim_2}
The bridge $X_{n_0}$ breaks up after finitely many pull-backs.
\end{claim}
\begin{proof}[Proof of Claim \ref{Claim_2}:]
Suppose that $X_n$ is a bridge from $\xi_1$ to $\xi_2$ for all $n\geq n_0$. Near the endpoints $\xi_1$ and $\xi_2$, $X_{n}$ consists of internal rays $R^{n}_1\subset U_{\xi_1}$ and $R^{n}_2\subset U_{\xi_2}$, respectively. Hence, we can express $X_n$ as $R^n_1\cup R^n_2\cup X'_n$, where $X'_n$ consists of a bounded number of edges of $\Delta_n$ by the previous considerations. Considering the situation in B\"ottcher coordinates as in the previous claim, we see that $R^n_1$ converges to a fixed ray $R^0_1\subset U_{\xi_1}$, and $R^n_2$ converges to a fixed ray $R^0_2\subset U_{\xi_2}$. Let $r_n$ be the non-fixed endpoint of $R^n_1$. Then, $r_n\to\infty$ as $n\to\infty$.

Since $\infty$ is a repelling fixed point with multiplier $\lambda>1$, there exists a branched covering map $\psi:(\C,0)\to (\Cc\setminus\{0\},\infty)$ such that $f(\psi(z))=\psi(\lambda z)$ for all $z\in\C$ (we may assume without loss of generality that $0\not\in \ol{V_0}$). Moreover, there exist a neighborhood $W$ of $\infty$ and a holomorphic branch $\phi$ of $\psi^{-1}$ on $W$ \cite[Corollary 8.10]{Milnor}. Let $N$ be sufficiently large so that  for all $n\geq N$, $r_n\in W$ and $\ol{X}'_n$ contains no critical values of $f$. Then, $\ol{X}'_{n+1}$ is the lift of $\ol{X}'_n$ under the branch of $f^{-1}$ that maps $r_n$ to $r_{n+1}$. Observe that a point $z\in\C$ is a critical point of $\psi$ if and only if $f$ has a critical point at $f^{\circ j}(\psi(z/\lambda^m))$ for some $j\leq m$, where $m$ is chosen large enough so that $z/\lambda^m\in \phi(W)$.  Therefore, none of the $X'_n$ for $n\geq N$ contain a critical value of $\psi$ and we can pull back $X'_n$ under the branch of $\psi^{-1}$ that maps $r_n$ to $\phi(r_n)$. Since $f$ commutes with multiplication by $\lambda$, this pull-back operation commutes with division by $\lambda$ and it follows that $X'_n$ converges uniformly to $\infty$ as a set as $n\to\infty$. Therefore, $X_n\to \ol{R}^0_1\cup\ol{R}^0_2$ uniformly as a set as $n\to\infty$. (Note that since $X_n$ is connected by assumption, we could have as well argued starting with the endpint of $R^n_2$.)

But this implies that if $n$ is sufficiently large, the arc $X_n$ separates the bounded set $C_1$ (and hence also $C_n\supset C_1$) from all accesses to $\infty$ of $V_n$. This contradicts Corollary \ref{Cor_FreeNest}.
\end{proof}

We can now finish the proof of Theorem~\ref{Thm_PolesConnect}.

If $X_{n_0}$ breaks up at time $n+1$, then consider its endpoint $\xi'_2\in f^{-1}(\{\xi_2\})$. If $\xi'_2\not\in \Delta_1$, then it must by assumption be in $C_1$. Thus, $X_{n+1}$ connects $C_{n}$ to $\Delta$. This means that $C_{n+1}\subset \Delta_{n+1}$ and we are done.

If $\xi'_2\in \Delta_1$, we can extend $X_{n+1}$ to a bridge by connecting $\xi'_2$ to a pole $p$ and on to another fixed point $\xi_3$ within $B_{n+1}$. Claim \ref{Claim_1} implies that $\xi_3\in\{\xi_1,\xi_2\}$, and it follows that $p\in X_n$, which must therefore be a multiple pole. But then, we can choose $\xi_3=\xi_2$ and have found a new bridge between $\xi_1$ and $\xi_2$. After pulling back this new bridge at most finitely many times, the case $\xi'_2\in\Delta_1$ cannot happen anymore and we arrive at a contradiction to our assumption.

If all $V_n$ surround several bounded components of $f^{-1}(\Delta)$, our arguments show that at least one of them is connected to $\Delta_k$ for some $k\in\N$. To finish the proof, it suffices to show that at time $k$, a new bridge is created that connects $\xi_1$ to some finite fixed point $\xi_3$. Then we can continue by induction.

To see that this bridge exists, let $C_1$ be the component of $f^{-1}(\Delta)$ that was connected to $\Delta_k$ by $X_k$. Observe that $C_1$ also contains a preimage of $\xi_1$. We can extend the arc $X_k$ from $\xi_1$ to  $\xi'_2\in C_1$ to this preimage of $\xi_1$ and then further to another preimage of $\xi_2$ etc., until we arrive at a finite fixed point, say $\xi_3$. As before, $\xi_3$ might equal $\xi_2$ if $X_k$ contains a multiple pole.

The above arguments apply to all components of $\Cc\setminus\Delta$ that surround a pole. 
\end{proof}

This finishes the proof of Theorem \ref{Thm_PreimagesConnected} under the finiteness condition \eqref{Prop_PCF}. In general, we can use a straightforward surgery construction to bring any Newton map into the desired form: all we require is finiteness of critical orbits in all basins of roots; there may well be infinite (or even periodic or preperiodic) critical orbits in the Julia set, or in attracting, parabolic, or Siegel components of the Fatou set. 

Within any immediate basin $U_\xi$, we may replace the attracting dynamics (which may involve several critical points converging to the root $\xi$) by dynamics modeled after $z\mapsto z^k$ within the unit disk. This surgery procedure does not affect the correctness of Theorem \ref{Thm_PreimagesConnected}. A similar procedure can assure that all critical points in the entire basin $B_\xi$ of $\xi$ land on the fixed point $\xi$ after finitely many steps. Details are standard and thus omitted; compare for example Shishikura~\cite{Shishikura}. This proves Theorem~\ref{Thm_PreimagesConnected} in the general case.

%%%%%%%%%%%%%%%%%%%%%%%%%%%
\section{The Newton Graph of a Newton Map}
\label{Sec_Graph}

In this section, we define abstract Newton graphs and use Theorem \ref{Thm_PreimagesConnected} to show that every postcritically fixed Newton map $f$ generates a unique abstract Newton graph in a natural way. We make some references to Thurston theory, which is discussed in more detail in Section \ref{Sec_Thurston}.

\subsection{Extending Maps on Finite Graphs}
The channel diagram motivates the definition of a Newton graph. For this, we first need to introduce some notation regarding maps on imbedded graphs and their extensions to $\S^2$, compare \cite[Chapter 6]{BFH}. We assume in the following without explicit mention that all graphs are imbedded into $\S^2$.

\begin{definition}[Regular Extension]
\label{Def_GraphMap}
Let $g:\Gamma_1\to\Gamma_2$ be a graph map. An orientation-preserving branched covering map $\ol{g}:\S^2\to\S^2$ is called a {\em regular extension} of $g$ if $\ol{g}|_{\Gamma_1}=g$ and $\ol{g}$ is injective on each component of $\S^2\setminus \Gamma_1$. 
\end{definition}

\begin{lemma}[Isotopic Graph Maps]
\label{Lem_IsotopGraphMaps}
{\cite[Corollary 6.3]{BFH}} Let $g,h:\Gamma_1\to\Gamma_2$ be two graph maps that coincide on the vertices of $\Gamma_1$ such that if $\gamma\subset\Gamma_1$ is an edge, then $g(\gamma)=h(\gamma)$ as a set. Suppose that $g$ and $h$ have regular extensions $\ol{g},\ol{h}:\S^2\to\S^2$. Then there exists a homeomorphism $\psi:\S^2\to\S^2$, isotopic to the identity relative the vertices of $\Gamma_1$, such that $\ol{g}=\ol{h}\circ \psi$.
\qed
\end{lemma}

Let $g:\Gamma_1\to\Gamma_2$ be a graph map. For the next proposition, we will assume without loss of generality that each vertex $v$ of $\Gamma_1$ has a neighborhood $U_v\subset\S^2$ such that all edges of $\Gamma_1$ that enter $U_v$ terminate at $V$; we may also assume that in local cordinates, $U_v$ is a round disk of radius $1$ centered at $v$, that all edges entering $U_v$ are radial lines and that $g|_{U_v}$ is length-preserving. We make analogous assumptions for $\Gamma_2$. Then, we can extend $g$ to each $U_v$ as in \cite{BFH}: for a vertex $v\in\Gamma_1$, let $\gamma_1$ and $\gamma_2$ be two adjacent edges ending there. In local coordinates, these are radial lines at angles, say, $\theta_1,\theta_2$ such that $0<\theta_2-\theta_1\leq 2\pi$ (if $v$ is an endpoint of $\Gamma_1$, then set $\theta_1=0$, $\theta_2=2\pi$). In the same way, choose arguments $\theta_1',\theta_2'$ for the image edges in $U_{g(v)}$ and extend $g$ to a map $\tilde{g}$ on $\Gamma_1\cup\bigcup_v U_v$ by setting 
\[
(\rho,\theta)\mapsto \left(\rho, \frac{\theta_2'-\theta_1'}{\theta_2-\theta_1}\cdot\theta\right)\;,
\]
where $(\rho,\theta)$ are polar coordinates in the sector bounded by the rays at $\theta_1$ and $\theta_2$. In other words, sectors are mapped onto sectors in an orientation-preserving way. Then, the following holds.

\begin{proposition}[Regular Extension]
\label{Prop_RegExt}
{\em \cite[Proposition 6.4]{BFH}} The map $g:\Gamma_1\to\Gamma_2$ has a regular extension 
if and only if for every vertex $y\in\Gamma_2$ and every component $U$ of $\S^2\setminus\Gamma_1$, the extension $\tilde{g}$ is injective on
\[
	\bigcup_{v\in g^{-1}(\{y\})} U_v \cap U\;.
\]
In this case, the regular extension $\ol{g}$ may have critical points only at the vertices of $\Gamma_1$.
\qed
\end{proposition}

%%%%%%%%%%%%%%%%%%%%%%%%%%%%%%%%%
\subsection{The Newton Graph}

With these preparations, we are ready to introduce the concept of an abstract Newton graph. It turns out that it carries enough information to uniquely characterize postcritically fixed Newton maps. 

\begin{definition}[Abstract Newton Graph]
\label{Def_NewtonGraph}
Let $\Gamma\subset\S^2$ be a connected graph, $\Gamma'$ the set of its vertices and $g:\Gamma\to\Gamma$ a graph map. The pair $(\Gamma,g)$ is called an {\em abstract Newton graph} if it satisfies the following conditions:
\begin{enumerate}
\item \label{Cond_Diagram} There exists $d_{\Gamma}\geq 3$ and an abstract channel diagram $\Delta\subsetneq\Gamma$ of degree $d_\Gamma$ such that $g$ fixes each vertex and each edge of $\Delta$.
\item  \label{Cond_Branch} If $v_0,\dots,v_{d_\Gamma}$ are the vertices of $\Delta$, then $v_i\in \ol{\Gamma\setminus\Delta}$ if and only if $i\neq 0$. Moreover, there are exactly $\deg_{v_i}(g)-1\geq 1$ edges in $\Delta$ that connect $v_i$ to $v_0$ for $i\neq 0$, where $\deg_x(g)$ denotes the local degree of $g$ at $x\in\Gamma'$.

\item \label{Cond_Degree} $\sum_{x\in\Gamma'} \left(\deg_x(g)-1\right) = 2d_{\Gamma}-2$.
\item \label{Cond_Enter} There exists $N_{\Gamma}\in\N$ such that $g^{\circ {N_\Gamma}}(\Gamma)\subset\Delta$, where $N_\Gamma$ is minimal such that $g^{\circ (N_{\Gamma}-1)}(x)\in\Delta$ for all $x\in\Gamma'$ with $\deg_x(g)>1$.
\item \label{Cond_Connected} The graph $\ol{\Gamma\setminus\Delta}$ is connected.
\item \label{Cond_Extension} For every vertex $y\in\Gamma'$ and every component $U$ of $\S^2\setminus\Gamma$, the extension $\tilde{g}$ is injective on 
\[
	\bigcup_{v\in g^{-1}(\{y\})} U_v \cap U\;.
\]
\item \label{Cond_Saturated} $\Gamma$ equals the component of $\ol{g}^{-N_{\Gamma}}(\Delta)$ that contains $\Delta$. 
\end{enumerate}
\end{definition}

If $(\Gamma,g)$ is an abstract Newton graph, $g$ can be extended to a branched covering map $\ol{g}:\S^2\to\S^2$ by (\ref{Cond_Extension}) and Proposition \ref{Prop_RegExt}. We use this implicitly in (\ref{Cond_Saturated}). Condition (\ref{Cond_Degree}) and the Riemann-Hurwitz formula ensure that $\ol{g}$ has degree $d_{\Gamma}$.
An immediate consequence of Lemma \ref{Lem_IsotopGraphMaps} is that $\ol{g}$ is unique up to Thurston equivalence.

We say that two abstract Newton graphs $(\Gamma_1,g_1)$ and $(\Gamma_2,g_2)$ are {\em equivalent} if there exists a graph homeomorphism $h:\Gamma_1\to\Gamma_2$ that preserves the cyclic order of edges at each vertex of $\Gamma_1$ and conjugates $g_1$ to $g_2$. 

Now we are ready to prove our first main result. Recall that for a Newton map $f$ with channel diagram $\Delta$, $\Delta_n$ denotes the component of $f^{-n}(\Delta)$ that contains $\Delta$. 

\begin{proof}[Proof of Theorem \ref{Thm_NewtonGraph}]
Let $\Delta$ be the channel diagram of $f$. First observe that $\Delta$ connects every fixed point of $f$ to $\infty$. Since $f$ is postcritically fixed, each critical point of $f$ is connected to some prepole by an iterated preimage of $\Delta$. Theorem \ref{Thm_PolesConnect} shows that there exists $n\in\N$ such that $\Delta_n$ contains all poles of $f$. Since $f$ is postcritically fixed, it follows by induction that there exists a minimal $N\in\N$ such that $\Delta_{N-1}$ contains all critical points of $f$.

It is easy to see that $(\Delta_N,f)$ satisfies all conditions of Definition \ref{Def_NewtonGraph} except possibly (\ref{Cond_Connected}). Note that if all critical points are in $\Delta_{N-1}$, we need to pull back one more step to ensure that condition (\ref{Cond_Degree}) is satisfied.

To show (\ref{Cond_Connected}), suppose by way of contradiction that the bounded set $\ol{\Delta_{N-1}\setminus\Delta}$ is not connected. Then, there exists an unbounded component $V$ of $\Cc\setminus\Delta_{N-1}$ that separates the plane, i.e.\ $V$ has at least two accesses to $\infty$. Let $W$ be a neighborhood of $\infty$ that is a round disk in linearizing coordinates and satisfies $W\cap \Delta_{N-1}\subset \Delta$. Let $V_1,\dots,V_k$ be the components of $V\cap W$. Then, $f$ acts injectively on each $V_i$ and there exists a branch $g_i$ of $f^{-1}$ that maps $V_i$ into itself (recall that $\infty$ is a repelling fixed point of $f$, so it is attracting for the $g_i$). By assumption, $V$ is simply connected and contains no critical values of $f$, so the $g_i$ extend to all of $V$ by holomorphic continuation on lines. 
Since $\Delta_{N-1}\subset\Delta_{N}$, we get $g_i(V)\subset V$ for all $i$. If there are $i\neq j$ such that $g_i(V)\cap g_j(V)\neq\emptyset$, then it follows that $g_i=g_j$ and we have found a holomorphic self-map of $V$ for which $\infty\in \partial V$ is attracting through two distinct accesses, contradicting the Denjoy-Wolff theorem \cite[Theorem 5.4]{Milnor}.

Hence, the $g_i(V)$ are pairwise disjoint and if $w\in \partial g_i(V)$ for some $i$, then $f(w)\in\partial V$, for otherwise the map $g_i$ would be defined in a neighborhood of $f(w)$. Hence $w\in f^{-1}(\Delta_{N-1})$ and since all $g_i(V)$ are open disks, we even get that $w\in \Delta_N$. It follows  that no component of $\Cc\setminus \Delta_N$ has more than one access to $\infty$, and $\ol{\Delta_N\setminus\Delta}$ is connected.

To prove the last claim, suppose that there exists a graph homeomorphism $h:\Delta_{1,N}\to \Delta_{2,N}$ such that $h(f_1(z)) = f_2(h(z))$ for all $z\in \Delta_{1,N}$. Since all complementary components of $\Delta_{1,N}$ are disks, we can extend $h$ to a homeomorphism $\ol{h}:\S^2\to\S^2$ that conjugates $f_1$ to $f_2$  up to isotopy relative $\Delta_{1,N}$ \cite[Lemma 6.1]{BFH}. By Theorem \ref{Thm_Thurston}, $f_1$ and $f_2$ are conjugate by a M\"obius transformation that fixes $\infty$.
\end{proof}

%%%%%%%%%%%%%%%%%%%%%%%%%%%%
\section{Abstract Newton Graphs Are Realized}
\label{Sec_Thurston}

In this section, we recall some fundamental notions of Thurston's characterization of rational maps. Then, we use Thurston's theorem to show that every abstract Newton graph is realized by a postcritically fixed Newton map, which is unique up to affine conjugation.

\subsection{Thurston's Criterion For Marked Branched Coverings}

Thurston's theorem provides a
necessary and sufficient condition on the existence of a rational
map with certain combinatorial behavior in terms of linear
maps generated from a (potentially very large) collection of simple closed curves.

\looseness -1
The notations and results in this section are based on \cite{DH} and \cite{PT}. 
Before we can state Thurston's criterion, we need several definitions. Recall that for a branched covering  map $g:\S^2\to\S^2$, $P_g$ denotes the postcritical set.

\begin{definition}[Marked Branched Covering]
A {\em marked branched covering} is a pair $(g,X)$, where $g:\S^2\to\S^2$ is a postcritically finite branched covering map and $X$ is a finite set containing $P_g$ such that $g(X)\subset X$.
\end{definition}
\begin{definition}[Thurston Equivalence]
Let $(g,X)$ and $(h,Y)$ be two marked branched coverings. We say that they are {\em (Thurston) equivalent} if there are two homeomorphisms $\phi_0,\phi_1:\S^2\to\S^2$ such that
\[	
	\phi_0\circ g = h\circ \phi_1 
\]
and there exists an isotopy $\Phi: [0,1]\times\S^2\to\S^2$ with $\Phi(0,.)=\phi_0$ and $\Phi(1,.)=\phi_1$ such that $\Phi(t,.)|_{X}$ is constant in $t\in [0,1]$ with $\Phi(t,X)=Y$.  
\end{definition}
If $(g,X)$ is a marked branched covering and $\gamma$ a simple closed curve in $\S^2\setminus X$, then the set $g^{-1}(\gamma)$ is a disjoint union of simple closed curves. 
\begin{definition}[Multicurve]
Let $(g,X)$ be a marked branched covering. We say that a simple closed curve $\gamma\subset\S^2$ is a {\em simple closed curve in $(\S^2,X)$} if $\gamma\subset\S^2\setminus X$. It is called {\em peripheral} if there exists a component of $\S^2\setminus\gamma$ that intersects $X$ in at most one point, and {\em non-peripheral} otherwise.

Two simple closed curves $\gamma_1,\gamma_2$ in $(\S^2,X)$ are called {\em isotopic (relative $X$)} (write $\gamma_1\simeq\gamma_2$) if there exists a continuous one-parameter family $\gamma_t$, $t\in [1,2]$, of such curves joining $\gamma_1$ to $\gamma_2$. We denote the isotopy class of $\gamma_1$ by $[\gamma_1]$.

A finite set $\Gamma=\{\gamma_1,\dots,\gamma_m\}$ of disjoint, non-peripheral and pairwise non-isotopic simple closed curves in $(\S^2,X)$ is called a {\em multicurve}. 
\end{definition}
\begin{definition}[Irreducible Thurston Obstruction]
Let $(g,X)$ be a marked bran\-ch\-ed covering and $\Gamma$ a multicurve. Denote by $\R^\Gamma$ the real vector space spanned by the isotopy classes of the curves in $\Gamma$. Then, we associate to $\Gamma$ its {\em Thurston transformation} $g_\Gamma:\R^\Gamma\to\R^\Gamma$ by specifying its action on representatives $\gamma\in\Gamma$ of basis elements:
\begin{equation}
\label{Eqn_ThurstonTransform}
	g_{\Gamma}(\gamma) := \sum_{\gamma'\subset g^{-1}(\gamma)} \frac{1}{\deg(g|_{\gamma'}:\gamma'\to\gamma)}[\gamma']\;.
\end{equation}
The sum is taken to be zero if there are no preimage components isotopic to a curve in $\Gamma$. 

The linear map given by equation (\ref{Eqn_ThurstonTransform}) is represented by a square matrix with non-negative entries and thus its largest eigenvalue $\lambda(\Gamma)$ is real and non-negative by the Perron-Frobenius theorem.

A square matrix $A_{i,j}\in\R^{n\times n}$ is called {\em irreducible} if for each $(i,j)$, there exists $k\geq 0$ such that $(A^k)_{i,j}>0$. We say that the multicurve $\Gamma$ is {\em irreducible} if the matrix representing $g_\Gamma$ is.

An irreducible multicurve $\Gamma$ is called an {\em irreducible (Thurston) obstruction} if $\lambda(\Gamma)\geq 1$.
\end{definition}

Now we are ready to state Thurston's theorem for marked branched coverings as given in \cite[Theorem 3.1]{PT} and proved in \cite{DH}. 
\begin{theorem}[Marked Thurston Theorem]
\label{Thm_Thurston}
Let $(g,X)$ be a marked branched covering with hyperbolic orbifold. It is Thurston equivalent to a marked rational map $(f,Y)$ if and only if it has no irreducible Thurston obstruction, i.e.\ if $\lambda(\Gamma)<1$  for each irreducible multicurve $\Gamma$. In this case, the rational map $f$ is unique up to automorphism of $\Cc$.
\qed
\end{theorem}

\begin{remark}
For the definition of a hyperbolic orbifold, see \cite{DH}. If $g$ has at least three fixed branch points, then it will have hyperbolic orbifold. In general, $\#P_g\geq 5$ suffices to make the orbifold of $g$ hyperbolic. 

Note that a marked rational map is in particular a
rational map, ``forgetting'' the marked set $Y$. 
\end{remark}
%%%%%%%%%%%%%%%%%%%%%%%%%%%%%%%%%%
\subsection{Arcs Intersecting Obstructions}

We present a theorem of Ke\-vin Pil\-grim and Tan Lei that is useful to show that certain marked branched coverings are equivalent to rational maps. Again, we first need to introduce some notation.

Let $(g,X)$ be a marked branched covering of degree $d\geq 3$.
\begin{definition}[Arc System]
An {\em arc} in $(\S^2,X)$ is a map $\alpha:[0,1]\to \S^2$ such
that $\alpha(\{0,1\})\subset X$, $\alpha((0,1))\cap X=\emptyset$
and $\alpha|_{(0,1)}$ is injective.
The notion of {\em isotopy relative $X$} extends to arcs and is
also denoted by $\simeq$.

A set of pairwise non-isotopic arcs in $(\S^2,X)$ is called an
{\em arc system}. Two arc systems $\Lambda,\Lambda'$ are {\em
isotopic} if each curve in $\Lambda$ is isotopic relative $X$ to a
unique element of $\Lambda'$ and vice versa.
\end{definition}
Note that arcs connect marked points (the endpoints of an arc need not be distinct) while
simple closed curves run around them. We will see that this leads to
intersection properties that will give us some control over the
location of possible Thurston obstructions. Since arcs and curves are only
defined up to isotopy, we make precise what we mean by arcs and
curves intersecting.
\begin{definition}[Intersection Number]
\label{Def_IntersectionNumber}
Let $\alpha$ and $\beta$ each be an arc or a
simple closed curve in $(\S^2,X)$. Their {\em intersection
number} is
\[
    \alpha\cdot\beta := \min_{\alpha'\simeq\alpha,
    \,\beta'\simeq\beta} \#\{(\alpha'\cap\beta')\setminus X\}\;.
\]
\end{definition}
\noindent
The intersection number extends bilinearly to arc systems and
multicurves.

If $\lambda$ is an arc in $(\S^2,X)$, then the closure of a
component of $g^{-1}(\lambda\setminus X)$ is called a {\em lift}
of $\lambda$. Each arc clearly has $d$ distinct lifts. If
$\Lambda$ is an arc system, an arc system $\tilde{\Lambda}$ is
called a {\em lift} of $\Lambda$ if each
$\tilde{\lambda}\in\tilde{\Lambda}$ is a lift of some
$\lambda\in\Lambda$.

If $\Lambda$ is an arc system, we introduce a linear map $g_{\Lambda}$ on the real vector space $\R^{\Lambda}$ similar as for multicurves: for $\lambda\in\Lambda$, set
\[
	g_{\Lambda}(\lambda):=\sum_{\lambda'\subset g^{-1}(\lambda)} [\lambda']\;,
\]
where $[\lambda']$ denotes the isotopy class of $\lambda'$ relative $X$. Again, the sum is taken to be zero if $\lambda$ has no preimages in the isotopy class of $\lambda'$. 
We say that $\Lambda$ is {\em irreducible} if the matrix representing $g_{\Lambda}$ is. 

Denote by $\tilde{\Lambda}(g^{\circ n})$ the union of those components of
$g^{-n}(\Lambda)$ that are isotopic to elements of $\Lambda$
relative $X$, and define $\tilde{\Gamma}(g^{\circ n})$ in an
analogous way. Note that if $\Lambda$ is irreducible, each
element of $\Lambda$ is isotopic to an element of
$\tilde{\Lambda}(g^{\circ n})$. 

The following theorem is Theorem 3.2 of \cite{PT}. It shows that up to isotopy, irreducible Thurston obstructions cannot intersect the preimages of irreducible arc systems (except possibly the arc systems themselves). We will use this theorem to show that the extended map of an abstract Newton graph is Thurston equivalent to a rational map.

\begin{theorem}[Arcs Intersecting Obstructions]
\label{Thm_ArcsInterOb} Let $(g,X)$ be a mar\-ked branched covering,
$\Gamma$ an irreducible Thurston obstruction and $\Lambda$ an
irreducible arc system. Suppose furthermore that
$\#(\Gamma\cap\Lambda)=\Gamma\cdot\Lambda$. Then, exactly one of
the following is true:
\begin{enumerate}
\item $\Gamma\cdot\Lambda=0$ and $\Gamma\cdot g^{-n}(\Lambda)=0$
for all $n\geq 1$.
\item $\Gamma\cdot\Lambda\neq 0$ and
for $n\geq 1$, each component of $\Gamma$ is isotopic to a unique
component of $\tilde{\Gamma}(g^{\circ n})$. The mapping $g^{\circ
n}:\tilde{\Gamma}(g^{\circ n})\to \Gamma$ is a homeomorphism and
$\tilde{\Gamma}(g^{\circ n})\cap (g^{-n}(\Lambda)\setminus
\tilde{\Lambda}(g^{\circ n}))=\emptyset$. The same is true when
interchanging the roles of $\Gamma$ and $\Lambda$.
\qed
\end{enumerate}
\end{theorem}

\subsection{The Realization of Abstract Newton Graphs}
We conclude with a proof of Theorem \ref{Thm_Realization}. If $\tilde{\Delta}$ is the channel diagram of the postcritically fixed Newton map $f$, recall that $\tilde{\Delta}_n$ denotes the component of $f^{-n}(\tilde{\Delta})$ that contains $\tilde{\Delta}$. By $\tilde{\Delta}'_n$ we denote the set of vertices of $\tilde{\Delta}_n$. 

\begin{proof}[Proof of Theorem \ref{Thm_Realization}]
Let $\Delta\subsetneq\Gamma$ be the abstract channel diagram of $\Gamma$.
First observe that by condition (\ref{Cond_Branch}) of Definition \ref{Def_NewtonGraph}, the vertices $v_1,\dots,v_{d_\Gamma}$ of $\Delta$ are branch points of $\ol{g}$. Since $d_{\Gamma}\geq 3$, $\ol{g}$ has hyperbolic orbifold and it suffices to show that $(\ol{g},X)$ has no irreducible Thurston obstruction: it then follows from Theorem \ref{Thm_Thurston} that $\ol{g}$ is Thurston equivalent to a rational map $f$ of degree $d_\Gamma$, which is unique up to M\"obius transformation. Then, $f$ has $d_\Gamma+1$ fixed points, $d_\Gamma$ of which are superattracting because $\ol{g}$ has marked fixed branch points $v_1,\dots,v_{d_\Gamma}$. The last fixed point is repelling \cite[Corollary 12.7 \& 14.5]{Milnor} and after possibly conjugating $f$ with a M\"obius transformation, we may assume that it is at $\infty$. Now it follows from Proposition \ref{Prop_Head} that $f$ is a Newton map. It is unique up to a M\"obius transformation fixing $\infty$, hence up to affine conjugacy.

So suppose by way of contradiction that $\Pi$ is an irreducible Thurston obstruction for $(\ol{g},X)$ and let $\gamma\in\Pi$. Then, $\gamma$ is a non-peripheral simple closed curve in $\S^2\setminus X$. 
It is easy to see that each edge $\lambda$ of $\Delta$ forms an irreducible arc system, hence Theorem \ref{Thm_ArcsInterOb} implies that $\gamma\cdot (\ol{g}^{-n}(\lambda)\setminus\lambda)=0$ for all $n\geq 1$. Since this is true for all edges of $\Delta$ and all vertices of $\Gamma$ are marked points, we get that $\gamma\cdot (\ol{\Gamma\setminus\Delta})=0$.
But since $\ol{\Gamma\setminus\Delta}$ is connected and contains $X\setminus\{v_0\}$, this means that $\gamma$ is peripheral, a contradiction.

In order to prove the last claim, note that $(\ol{g}_1,\Gamma'_1)$ and $(\ol{g}_2,\Gamma'_2)$ are Thurston equivalent as marked branched coverings. Let $h:(\S^2, \Gamma'_1)\to (\S^2,\Gamma'_2)$ be a homeomorphism that conjugates $\ol{g}_1$ to $\ol{g}_2$ on $\Gamma'_1$. If $e$ is an edge of $\Gamma_1$ with endpoints $x_1,x_2\in \Gamma'_1$, then $h(e)$ connects $h(x_1)$ with $h(x_2)$. Moreover, $h$ preserves the cyclic order at each vertex of $\Gamma_1$, because it is a homeomorphism of $\S^2$. So if $h':h(\Gamma_1)\to\Gamma_2$ is a homeomorphism that maps each $h(e)$ to the edge of $\Gamma_2$ between $h(x_1)$ and $h(x_2)$, then $h'\circ h$ realizes an equivalence between the two abstract Newton graphs. 
\end{proof}

\section{Acknowledgements}
I would like to thank the Fields Institute for Research in Mathematical Science in Toronto, Canada, for its hospitality while this paper evolved. I thank Dierk Schleicher for his everlasting support and for many fruitful discussions about the combinatorics of Newton maps. I also thank Tan Lei for her comments that helped to improve this paper.

%%%%%%%%%%%%%%%%%%%%%%%%%%%


\begin{thebibliography}{BFH}
\bibitem[BFH]{BFH} B. Bielefeld, Y. Fisher \and J. Hubbard: {\it The classification of critically preperiodic polynomials as dynamical systems}. J. Amer. Math. Soc. {\bf 5} (4) (1992), 721--762.
\bibitem[DH] {DH} A. Douady \and J. Hubbard: {\it A proof of Thurston's topological characterization
of rational functions.} Acta Math. {\bf 171} (1993), 263--297.
\bibitem[He] {Head} J. Head: {\it The combinatorics of Newton's method
for cubic polynomials}, Thesis Cornell University (1987).
\bibitem[HSS]{HSS} J. Hubbard, D. Schleicher \and S. Sutherland: {\it How to find all roots of complex polynomials by Newton's method.} Invent. Math. {\bf 146} (2001), 1--33.
\bibitem[Lu]{Luo} J. Luo: {\it Newton's method for polynomials with one inflection value}, preprint, Cornell University 1993.
\bibitem[Mi1]{Milnor2} J. Milnor, {\it Geometry and dynamics of quadratic rational maps}, with an appendix by the author and Tan Lei. Exp. Math. 2 (1993), no.\ 1, 37--83.
\bibitem[Mi2] {Milnor} J. Milnor: {\it Dynamics in One Complex
Variable}, Vieweg (2000).
\bibitem[PT] {PT} K. Pilgrim \and T. Lei: {\it Combining rational maps and controlling
obstructions.} Ergodic Theory Dynam. Systems {\bf 18} (1998)
221--245.
\bibitem[Pr] {Przytycki} F. Przytycki: {\it Remarks on the simple connectedness
of basins of sinks for iterations of rational maps},
preprint, Polish Academy of Sciences, Warsaw (1987).
\bibitem[Ro]{Roesch} P. Roesch, {\it Topologie locale des m\'etodes de Newtopn cubiques: plan dynamique.} C. R. Acad. Sci. Paris S\'erie I 326 (1998), 1221--1226.
\bibitem[RS] {RS} J. R\"uckert \and D. Schleicher, {\it On Newton's method for entire functions}, J. London Math. Soc., to appear. ArXiv:math.DS/0505652.
\bibitem[Sh]{Shishikura} M. Shishikura: {\it The connectivity of the Julia set and fixed points}, preprint, IHES (1990).
\bibitem[TL]{TanLei} Tan Lei: {\it Branched coverings and cubic Newton maps}. Fund. Math. {\bf 154} (1997), 207--260.
\end{thebibliography}
\end{document}